\documentclass[leqno,12pt]{article}
\usepackage{amssymb,amsmath}

\newtheorem{theorem}{Theorem}
\newtheorem{proposition}[theorem]{Proposition}
\newtheorem{lemma}[theorem]{Lemma}
\newtheorem{definition}[theorem]{Definition}
\newtheorem{corollary}[theorem]{Corollary}
\newenvironment{proof}{\noindent{\bf Proof. }}{\par}

\renewcommand{\Box}
        {\hbox{\hskip 1pt \vrule width4pt height 6pt depth 1.5pt \hskip 1pt}}
\newcommand{\qed}{\relax{\ifhmode\unskip\nobreak$\kern10pt\Box$\fi
        \ifmmodeh\ifinner\else\hskip5pt\fi \kern10pt\Box\fi}\relax}

\newcommand{\C}{{\mathbb C}}
\newcommand{\R}{{\mathbb R}}
\newcommand{\Z}{{\mathbb Z}}

\newcommand{\CB}{{\cal B}}

\newcommand{\CH}{{\cal H}}

\newcommand{\CO}{{\cal O}}

\title{Arrangement of hyperplanes II: \\
Szenes formula and Eisenstein series}
\author{Michel Brion and Mich{\`e}le Vergne}

\date{}
\begin{document}
\maketitle

{\bf Dedicated to Victor Guillemin, for 
his 60th birthday}

\bigskip

{\bf Abstract:}
The aim of this article is to generalize in several variables
some formulae for Eisenstein series in one variable.
In particular, we relate Szenes formula to Eisenstein series and 
we give another proof of it.

\section{Introduction}
Consider a sequence $(\alpha_1,\alpha_2,\ldots,\alpha_k)$
of linear forms in $r$ complex variables, with integral coefficients. 
The linear forms $\alpha_j$ need not be distinct.
For example $r=2$ and $\alpha_1=\alpha_2=z_1$,
$\alpha_3=\alpha_4=z_2$, $\alpha_5=\alpha_6=z_1+z_2$.
For any such sequence, Zagier \cite{Z} introduced the series
$$
\sum_{n\in \Z^r, \langle \alpha_j,n\rangle \neq 0}
\frac{1}{\prod_{j=1}^k\langle \alpha_j,n\rangle }.
$$
Assuming convergence, its sum is a rational multiple of $\pi^k$. 
For example \cite{Z}, we have
$$
\sum_{n_1\neq 0, n_2\neq 0, n_1+n_2\neq 0}
\frac{1}{n_1^2n_2^2 (n_1+n_2)^2}=\frac{(2\pi)^6}{30240}.
$$
These numbers are natural multidimensional generalizations of the
value of the Riemann zeta function at even integers.
A.~Szenes gave in (\cite{S}, Theorem 4.4) a residue formula for 
these numbers, relating them to Bernoulli numbers. The formula of
Szenes \cite{S} is the multidimensional analogue of the residue
formula 
$$
\sum_{n\neq 0} \frac{1}{n^{2l}}=(2\pi)^{2l}\frac{B_{2l}}{(2l)!}
=(-1)^l (2\pi)^{2l} Res_{z=0}\left(\frac{1}{z^{2l}(1-e^z)}\right).
$$

A motivation for computing such sums comes from the work of E.~Witten
\cite{W}: In the special case where $\alpha_j$ are the positive roots
of a compact connected Lie group $G$, each of these roots being
repeated with multiplicity $2g-2$, Witten expressed the symplectic
volume of the space of homomorphisms of the fundamental group of a
Riemann surface of genus $g$ into $G$, in terms of these sums.
In \cite{JK}, L.~Jeffrey and F.~Kirwan proved a special case of Szenes
formula leading to the explicit computation of this symplectic volume,
when $G$ is SU($n$).

Our interest in such series comes from a different motivation. Let us
consider first the one-dimensional case. By Poisson formula, the
convergent series, for $Re(z)>0$, 
$\sum_{m=1}^{\infty} me^{-mz}$ is also equal to
$\sum_{n\in\Z}1/(z+2i\pi n)^2$. Similarly, sums of products of
polynomial functions with exponential functions over all integral
points of a $r$-dimensional rational convex cone are related to
functions of $r$ complex variables of the form
$$
\psi(z)=\sum_{n\in\Z^r} 
\frac{1}{\prod_{j=1}^k \langle\alpha_j,z+2i\pi n\rangle}.
$$
When this series is not convergent, introduce the oscillating factor 
$e^{\langle t,2i\pi n\rangle}$ and define the Eisenstein series
$$
\psi(t,z)=\sum_{n\in\Z^r}
\frac{e^{\langle t,z+2i\pi n\rangle}}
{\prod_{j=1}^k \langle \alpha_j,z+2i\pi n\rangle},
$$
a generalized function of $t\in\R^n$. 

In Section 3, we construct a decomposition of an open dense subset of
$\R^n$ into alcoves such that $t\mapsto \psi(t,z)$ is given on each
alcove by a polynomial in $t$, with coefficients rational functions of
$e^z$. Our first theorem (Theorem \ref{ke}) gives an explicit residue
formula for $\psi(t,z)$. It follows easily from the obvious behaviour
of $\psi(t,z)$ under differentiation in $z$.
 
This formula allows to give a residual meaning ``$\psi(t,0)$'' for the
value of $\psi(t,z)$ at $z=0$, although $\psi(t,z)$ has clearly poles
along all hyperplanes $\langle \alpha_j,z\rangle=0$. An alternative
way to define $\psi(t,0)$ is to remove all infinities $1/\alpha_j$ in
the series
$$
\psi(t,0)=\sum_{n\in\Z^r}
\frac{e^{\langle t,2i\pi n\rangle}}
{\prod_{j=1}^k \langle\alpha_j,2i\pi n\rangle}.
$$ 
Indeed, we prove that the residue formula for ``$\psi(t,0)$''
coincides with the renormalized sum~:
$$
``\psi(t,0)''=\sum_{n\in\Z^r,\langle\alpha_j,n\rangle \neq 0}
\frac{e^{\langle t,2i\pi n\rangle}}
{\prod_{j=1}^k \langle \alpha_j,2i\pi n\rangle}.
$$
This equality gives another proof of Szenes residue formula, as a
``limit'' of a natural formula for $\psi(t,z)$ when $z\to 0$ along a
generic line.

To illustrate our method, let us consider the one-dimensional case.
For $k\geq 2$, we can define the Eisenstein series 
$$
E_k(z)=\sum_{n\in \Z}\frac{1}{(z+2i\pi n)^k}.
$$
Clearly, $E_k(z)$ is periodic in $z$ with respect to translation 
by the lattice $2i\pi\Z$. From the residue theorem, when $y$ is not in
$2i\pi \Z$, we have the kernel formula: 
\begin{equation}\label{ker}
E_k(y)=Res_{z=0} \left(\frac{1}{z^k(1-e^{z-y})}\right).
\end{equation}
Observe that the right-hand side has a meaning when $y=0$, and equals
by definition the Bernoulli number $B_k/k!$. The function 
$$
E_k(y)=\frac{1}{y^k}+\sum_{n\in \Z, n\neq 0}\frac{1}{(y+2i\pi n)^k}
$$ 
has a Laurent expansion at $y=0$, with $1/y^k$ as Laurent negative
part. We see from the residue formula that the constant term 
$CT(E_k)=\sum_{n\in\Z,n\neq 0} \frac{1}{(2i\pi n)^k}$
equals
$Res_{z=0} \left(\frac{1}{z^k(1-e^z)}\right)$.

In view of this example, we call the value ``$\psi(t,0)$''
of $\psi(t,y)$ at $y=0$ the constant term of the Eisenstein series
$$
\sum_{n\in\Z^r} \frac{e^{\langle t,z+2i\pi n\rangle}}
{\prod_{j=1}^k \langle \alpha_j,z+2i\pi n\rangle}.
$$

We thank A.~Szenes and the referees of our paper for several
suggestions. 

\section{Kernel formula}

In this section, we recall briefly results of 
\cite{BV} with slightly modified notation. Let $V$ be a 
$r$-dimensional complex vector space.
Let $V^*$ be the dual vector space and let
$\Delta\subset V^*$ be a finite subset of non-zero linear forms.
Each $\alpha\in \Delta$ determines a hyperplane 
$\{\alpha=0\}$ in $V$. Consider the hyperplane arrangement
$$
\CH=\bigcup_{\alpha\in \Delta}\{\alpha=0\}.
$$
An element $z\in V$ is called {\bf regular}
if $z$ is not in $\CH$. If $S$ is a subset of $V$, we write
$S_{reg}$ for the set of regular elements in $S$.
The ring $R_\Delta$ of rational functions with poles on $\CH$
is the ring 
$\Delta^{-1}S(V^*)$ 
generated by the ring 
$S(V^*)$ of polynomial functions on $V$, together with
inverses of the linear functions $\alpha\in \Delta$.
The ring $R_{\Delta}$ has a $\Z$-gradation by the homogeneous degree 
which can be positive or negative.
Elements of $R_{\Delta}$ are defined on the open subset $V_{reg}$.
(Our notation differs from \cite{BV} in that the roles of $V$ and
$V^*$ are interchanged.)

In the one variable case, the function $1/z$
is the unique function which cannot be obtained as a derivative.
There is a similar description 
of a complement space to the space of derivatives in the ring
$R_{\Delta}$ that we recall now.

A subset $\sigma$ of $\Delta$ is called a {\bf basis of} $\Delta$, if
the elements $\alpha\in \sigma$ form a basis of $V$. We denote by
$\CB(\Delta)$ the set of bases of $\Delta$. An {\bf ordered basis} is
a sequence $(\alpha_1,\alpha_2,\ldots,\alpha_r)$
of elements of $\Delta$  such that the underlying set 
is a basis. We denote by $O\CB(\Delta)$ the set of ordered bases.

For $\sigma\in \CB(\Delta)$, set 
$$
\phi_{\sigma}(z):=\frac{1}{\prod_{\alpha\in \sigma}\alpha(z)}.
$$
We call $\phi_{\sigma}$ a {\bf simple fraction}. Setting 
$z_j=\langle z,\alpha_j\rangle $, we have
$$
\phi_{\sigma}(z)=\frac{1}{z_1 z_2\cdots z_r}.
$$
\begin{definition}
The subspace $S_{\Delta}$ of $R_{\Delta}$ spanned by the elements 
$\phi_{\sigma},~ \sigma\in\CB(\Delta),$ will be called the space of
{\bf simple elements} of $R_{\Delta}$:
$$
S_{\Delta}=\sum_{\sigma\in \CB(\Delta)}\C \phi_{\sigma}.
$$
\end{definition}

The space $S_{\Delta}$ consists of homogeneous rational functions of 
degree $-r$. However, not every homogeneous element of degree 
$-r$ of $R_{\Delta}$ is in $S_{\Delta}$ (for example, in the
preceding notation if $r\geq 2$, both functions $1/z_1^r$ and 
$z_2/z_1^{r+1}$ are not in $S_\Delta$). 
Furthermore  we must be careful, as the elements 
$\phi_{\sigma}$ may be linearly dependent:
for example, if $V=\C^2$ and 
$\Delta=\{z_1, z_2, z_1+z_2\}$,
we have 
$$
S_\Delta=
\C\frac{1}{z_1 z_2}+\C\frac{1}{z_1(z_1+z_2)}+\C\frac{1}{z_2(z_1+z_2)}
$$
and we have the relation
$$
\frac{1}{z_1 z_2}=\frac{1}{z_1(z_1+z_2)}+\frac{1}{z_2(z_1+z_2)}.
$$

A description due to Orlik and Solomon of all linear relations between
the elements $\phi_{\sigma}$ is given in \cite{BV} Proposition 13.

\begin{definition}
A {\bf basis} $B$ of $\CB(\Delta)$ is a subset of $\CB(\Delta)$ such
that the elements ($\phi_{\sigma}$, $\sigma\in B$), form a basis of
$S_{\Delta}$:
$$
S_{\Delta}=\bigoplus_{\sigma\in B}\C \phi_{\sigma}.
$$
\end{definition}

We let elements $v$ of $V$ act on $R_{\Delta}$ by differentiation:
$$
(\partial(v)f)(z):=
\frac{d}{d\epsilon}f(z+\epsilon v)|_{\epsilon=0}.
$$
Then the following holds (\cite{BV}, Proposition 7.)
\begin{theorem}
$$
R_\Delta= \partial(V) R_{\Delta}\oplus S_{\Delta}.
$$
\end{theorem}
Thus we see that only simple fractions cannot be obtained as derivatives.

As a corollary of this decomposition, 
we can define the projection map
$$
Res_{\Delta}: R_{\Delta}\to S_{\Delta}.
$$
The projection $Res_{\Delta}f(z)$
of a function $f(z)$ is a function of $z$ that we  called the
{\bf Jeffrey-Kirwan residue} of $f$. By definition, this function 
can be expressed as a linear combination 
of the simple fractions $\phi_{\sigma}$.
The main property of the map $Res_{\Delta}$ is that it
vanishes on derivatives, so that
for $v\in V$, $f, g\in R_{\Delta}$:
\begin{equation}\label{commute}
Res_\Delta\left((\partial(v)f)\, g\right)=
-Res_{\Delta}\left(f \, (\partial(v)g)\right).
\end{equation}

If $o\sigma\in O\CB(\Delta)$ is an ordered basis, an important functional 
$Res^{o\sigma}$ can be defined on  $R_{\Delta}$:
the {\bf iterated residue} with respect to the ordered basis 
$o\sigma$. If we write an element 
$z\in V$ on the basis $o\sigma=(\alpha_1,\alpha_2,\ldots,\alpha_r)$
as $z=(z_1,\ldots,z_r)$, then
$$
Res^{o\sigma}(f)= 
Res_{z_1=0}(Res_{z_2=0}\ldots(Res_{z_r=0}f(z_1,z_2,\ldots,z_r))\ldots).
$$
The map $Res^{o\sigma}$ depends on the order $o\sigma$ chosen on
$\sigma$ and not only on the basis $\sigma$ underlying $o\sigma$. The
restriction of the functional $Res^{o\sigma}$ to $S_{\Delta}$ is
called $r^{o\sigma}$. 
We have 
\begin{equation}\label{compatible}
Res^{o\sigma}=r^{o\sigma}Res_{\Delta}.
\end{equation}
Indeed, we have only to check that $Res^{o\sigma}$ vanishes on 
derivatives. If $o\sigma=(\alpha_1,\alpha_2,\ldots,\alpha_r)$ and
$z=(z_1,\ldots,z_r)$, the iterated residue $Res^{o\sigma}$ will vanish
at the step $Res_{z_j=0}$ on $\frac{\partial}{\partial z_j}R_{\Delta}$.

Recall the following definition of A.~Szenes (\cite{S}, Definition
3.3).

\begin{definition}
A {\bf diagonal basis} is a subset $OB$ of $O\CB(\Delta)$
such that 

1) The set of underlying (unordered) bases forms a basis $B$ of
$\CB(\Delta)$.

2) The dual basis to the basis 
($\phi_{\sigma}$, $o\sigma\in OB$) is the set of linear forms 
($r^{o\sigma}$, $o\sigma\in OB$):
$$
r^{o\tau}(\phi_{\sigma})=\delta_{\sigma}^{\tau}.
$$
\end{definition}

In Proposition 3.4 of \cite{S}, it is proved that a total order on
$\Delta$ gives rise to a diagonal basis (this is reproved in more
detail in \cite{BV}, Proposition 14.) 
\bigskip

In the one dimensional case, $S_{\Delta}=\C z^{-1}$, and
the space $G=\sum_{k\leq -1} \C z^{_k}$ of negative Laurent series 
is the space obtained from the function $1/z$
by successive derivations. In the case of several variables, 
we can also characterize the space generated by simple fractions under
differentiation. 

Let $\kappa$ be a sequence of (not necessarily distinct)
elements of $\Delta$. The sequence $\kappa$ is called 
{\bf generating} if the $\alpha\in \kappa$ generate the vector space $V^*$.

We denote by $G_{\Delta}$ the subspace of $R_\Delta$ spanned by the 
$$
\phi_{\kappa}:=\frac{1}{\prod_{\alpha\in \kappa}\alpha}
$$
where $\kappa$ is a generating sequence. Finally, we denote by $S(V)$ 
the ring of differential operators on $V$, with constant
coefficients. This ring acts on $S(V^*)$ and on $R_{\Delta}$.

\begin{proposition} (\cite{BV}, Theorem 1)
The space $G_{\Delta}$ is the $S(V)$-submodule of $R_{\Delta}$
generated by $S_{\Delta}$.
\end{proposition}
For example, if $\Delta=\{z_1, z_2, z_1+z_2\}$, we have  
$$
\frac{1}{z_1 z_2 (z_1+z_2)}=
-\frac{\partial}{\partial z_1}(\frac{1}{z_1 z_2}) +
(\frac{\partial}{\partial z_1}-\frac{\partial}{\partial z_2})
(\frac{1}{z_1(z_1+z_2)}) .
$$

In particular, every element of $G_{\Delta}$
can be expressed as a linear combination of elements 
$$
\frac{1}{\prod_{\alpha\in \sigma}\alpha^{n_\alpha}}
$$
where $\sigma$ is a basis and the $n_\alpha$ are positive integers.

For example, the above equality is equivalent to
$$
\frac{1}{z_1 z_2 (z_1+z_2)}=
\frac{1}{z_1^2 z_2} - \frac{1}{z_1^2(z_1+z_2)}.
$$ 

\bigskip

The ring $S(V^*)$ operates by multiplication on $R_{\Delta}$.
It is also useful to consider the action of the ring 
${\cal D}(V)$ of differential operators with polynomial coefficients,
generated by $S(V)$ and $S(V^*)$. The following lemma is an obvious
corollary of the description of $G_\Delta$.

\begin{lemma}\label{mod}
The space $R_\Delta$ is generated by $G_\Delta$
as an $S(V^*)$-module. It is generated by 
$S_\Delta$ as a ${\cal D}(V)$-module.
\end{lemma}

Consider now the space $\CO$
of holomorphic functions on $V$ defined in 
a neighborhood of $0$. Let $\CO_{\Delta}=\Delta^{-1} \CO$
be the space of meromorphic functions in a neighborhood 
of $0$, with denominators products of elements of $\Delta$. 
The space $\CO_{\Delta}$ is a module for the action of differential
operators with constant coefficients. 
Via the Taylor series at the origin of elements 
of $\CO$, the residue $Res_\Delta f(z)$ has still a meaning if 
$f(z)\in \CO_{\Delta}$; indeed, $Res_{\Delta}f(z)=0$ if 
$f\in R_{\Delta}$ is homogeneous of degree $\neq -r$.
 
If $y\in V$ is sufficiently near $0$ and $f\in \CO_\Delta$, the
function 
$$
({\cal T}(y)f)(z):= f(z-y)
$$
is still an element of $\CO_{\Delta}$.
If moreover $y$ is regular, then $f(z-y)$ is defined for $z=0$, and
thus is an element of $\CO$. 
 
If $f\in R_{\Delta}$, we denote by $m(f)$ the operator of
multiplication by $f$: 
$$
(m(f)\phi)(z):=f(z)\phi(z).
$$
It operates on $\CO_{\Delta}$. Finally, we denote by $C$ the operator 
$$
(C f)(z)~:=f(-z)
$$
on $\CO_{\Delta}$.

\begin{theorem}\label{kern} (Kernel theorem). Let 
$A:R_{\Delta}\to \CO_{\Delta}$
be an operator commuting with the action 
of differential operators with constant coefficients.
For $y\in V$ regular, sufficiently 
near $0$ and for $f\in G_{\Delta}$, we have the formula
$$
(Af)(y) =Tr_{S_{\Delta}}
\left(Res_{\Delta} m(f)C{\cal T}(y) A Res_{\Delta}\right).
$$
More explicitly, choose a basis $B$ of $\CB(\Delta)$ and let
($\phi^{\sigma}, \sigma\in B$) be the basis of $S_{\Delta}^*$ dual to
the basis ($\phi_{\sigma},\sigma\in B$) of $S_{\Delta}$. Then we have
the kernel formula:
$$
(Af)(y)=\sum_{\sigma\in B} \langle \phi^{\sigma},
Res_{\Delta}\left(f(z)A_{\sigma}(y-z)\right)\rangle 
$$
where $A_{\sigma}(z)=A(\phi_{\sigma})(z)$.
\end{theorem}
 
Concretely, this formula means the following. Let $f$ be
homogeneous of degree $d$.  We fix $y$ regular and small.
The function $z\mapsto A_{\sigma}(y-z)$ is defined near $z=0$.
The Jeffrey-Kirwan residue $Res_{\Delta}$
of the function $z\mapsto f(z)A_{\sigma}(y-z)$
is a function of $z$ belonging to the space $S_\Delta$.
We pair it with the linear form $\phi^{\sigma}$ on $S_\Delta$
and obtain a certain complex number depending on $y$.
More precisely, consider the Taylor expansion
$$
A_{\sigma}(y-z)=A_{\sigma}(y)+\sum_{j=1}^{\infty}A_{\sigma}^j(y,z) 
$$
where $A_{\sigma}^j(y,z)$ is the part of the Taylor expansion
at $0$ of the holomorphic function $z\mapsto A_{\sigma}(y-z)$,
which is homogeneous of degree $j$ in $z$. We have 
$$
A_{\sigma}^j(y,z)=(-1)^j\sum_{(k), |(k)|=j}
A_{\sigma}^{(k)}(y) \frac{z^{(k)}}{(k)!}
$$ 
where $(k)=(k_1,\ldots,k_r)$ is a multi-index, and 
$A_{\sigma}^{(k)}(y) =
((\frac{\partial}{\partial y})^{(k)}A_{\sigma})(y)$.
Then, as the Jeffrey-Kirwan residue vanishes on homogeneous terms 
of degree not equal to $-r$, we obtain
$$\displaylines{
Res_{\Delta}\left(f(z)A_{\sigma}(y-z)\right)= 
Res_{\Delta}\left(f(z)A_{\sigma}^{-d-r}(y,z)\right)
\hfill\cr\hfill
=(-1)^{d+r}\sum_{(k), |(k)|=-d-r}
A_{\sigma}^{(k)}(y) Res_{\Delta}\left(f(z)\frac{z^{(k)}}{(k)!}\right).
\cr}
$$ 
Thus, 
$\langle \phi^{\sigma},
Res_{\Delta}\left(f(z)A_{\sigma}(y,z)\right)\rangle$
is equal to
$$
(-1)^{d+r}\sum_{(k), |(k)|=-d-r} 
A_{\sigma}^{(k)}(y) \langle\phi^{\sigma},
Res_{\Delta}\left(f(z)\frac{z^{(k)}}{(k)!}\right)\rangle.
$$ 
Set $c_{\sigma}^{(k)}(f)=\langle \phi^{\sigma},
Res_{\Delta}\left(f(z)\frac{z^{(k)}}{(k)!}\right)\rangle $.
Let $P_{\sigma}^f(\frac{\partial}{\partial y})$ 
be the differential operator with constant coefficients defined by
$$
P_{\sigma}^f(\frac{\partial}{\partial y})=(-1)^{d+r}
\sum_{(k), |(k)|=-d-r}c_{\sigma}^{(k)}(f)
(\frac{\partial}{\partial y})^{(k)}.
$$
Then $P_{\sigma}^f$ depends linearly on $f$, and 
$$
\langle \phi^{\sigma},
Res_{\Delta}\left(f(z)A_{\sigma}(y-z)\right)\rangle 
=(P_{\sigma}^f(\frac{\partial}{\partial y}) A_{\sigma})(y).
$$
The claim of the theorem is that 
$$
(Af)(y)=\sum_{\sigma\in B} 
P_{\sigma}^f(\frac{\partial}{\partial y})\cdot A_{\sigma}(y).
$$
We now prove this theorem.

\begin{proof}
Define an operator $A':R_{\Delta}\to\CO_{\Delta}$ by
$$
(A'f)(y)=
\sum_{\sigma\in B}\langle \phi^{\sigma},
Res_{\Delta}\left(f(z)A_{\sigma}(y-z)\right)\rangle.
$$
We first check that $A'$ commutes with the action 
of differential operators with constant coefficients.
Using the equation 
$$
(\partial_y(v)\phi)(y-z)=-(\partial_z(v)\phi)(y-z)
$$
and the main property (\ref{commute}) of $Res_{\Delta}$, we obtain 
$$
\displaylines{
\partial_y(v)\cdot\langle \phi^{\sigma},
Res_{\Delta}\left(f(z)A_{\sigma}(y-z)\right)\rangle 
=\langle \phi^{\sigma},Res_{\Delta}
\left(f(z)(\partial_y(v)\cdot A_{\sigma}(y-z))\right)\rangle
\hfill\cr\hfill
=-\langle \phi^{\sigma},Res_{\Delta}
\left(f(z)(\partial_z(v)\cdot A_{\sigma}(y-z))\right)\rangle 
=\langle \phi^{\sigma},Res_{\Delta}
\left(\partial_z(v)\cdot f) A_{\sigma}(y-z)\right)\rangle.
\cr}$$

It remains to see that $A$ and $A'$ coincide on $S_{\Delta}$. For
this, we will use the following formula. If $P$ is a polynomial 
and $\phi$ a simple fraction, then
\begin{equation}\label{0}
Res_\Delta(P \phi)=P(0)\phi.
\end{equation}
To see this, recall that the function $\phi$ is homogeneous 
of degree $-r$. As $P\in S(V^*)$, $P-P(0)$ is a sum of homogeneous
terms of positive degree. Thus, for homogeneity 
reasons, $Res_{\Delta}((P-P(0))\phi)=0$.

Let $y$ be regular and let $\sigma,\tau\in B$. As the function 
$z\to A_{\sigma}(y-z)$ is an element of $\CO$,
we obtain by formula (\ref{0}),
$$
Res_{\Delta}\left(\phi_{\tau}(z)A_{\sigma}(y-z)\right)
=A_{\sigma}(y) \phi_{\tau}(z).
$$
Thus 
$$\displaylines{
A'(\phi_{\tau})(y)=\sum_{\sigma\in B}
\langle \phi^{\sigma},Res_{\Delta}
\left(\phi_{\tau}(z)A_{\sigma}(y-z)\right)\rangle
\hfill\cr\hfill
=\sum_{\sigma\in B}\langle\phi^{\sigma},
\phi_{\tau}\rangle A_{\sigma}(y)
=\sum_{\sigma\in B}\delta_{\sigma}^{\tau}A_{\sigma}(y)=
A_{\tau}(y)= A(\phi_{\tau})(y).\cr}
$$
\end{proof}

Choosing a diagonal basis $OB$ and using Equation \ref{compatible}, we
obtain an iterated residue fomula for $(Af)(y)$:

\begin{corollary}\label{iterated}
For any diagonal basis $OB$ of $\CB(\Delta)$, we have for 
$f\in G_{\Delta}$:
$$
(Af)(y)=\sum_{o\sigma\in OB} Res^{o\sigma}(f(z)A_{\sigma}(y-z))
$$ 
where $A_{\sigma}(z)=A(\phi_{\sigma})(z)$.
\end{corollary}

Corollary \ref{iterated} applies to the identity operator 
$A:R_{\Delta}\to R_{\Delta}$. If $f\in G_{\Delta}$, we obtain 
$f(y)=\sum_{o\sigma\in OB} Res^{o\sigma}(f(z)\phi_{\sigma}(y-z))$. 
But if $f\in NG_{\Delta}$ then clearly
$Res^{o\sigma}(f(z)\phi_{\sigma}(y-z))=0$ as the Taylor series of
$f(z)\phi_{\sigma}(y-z)$ at $z=0$ is also in $NG_{\Delta}$. As a
consequence, we obtain a formula for the Jeffrey-Kirwan residue as a 
function of iterated residues:

\begin{lemma}
For any $f\in R_{\Delta}$, we have
$$
(Res_{\Delta}f)(y)=\sum_{o\sigma\in OB}
Res^{o\sigma}(f)\phi_{\sigma}(y).
$$ 
\end{lemma}

\bigskip
Similarly, if $Z:R_{\Delta}\to \CO$ is an operator commuting 
with the action of differential operators with constant coefficients,
the formula 
$$
Z(f)(y)=Tr_{S_{\Delta}}
\left(Res_{\Delta}m(f)C{\cal T}(y)Z Res_{\Delta}\right)
$$
is valid for {\bf all} elements $y\in V$ sufficiently near $0$
and for all $f\in G_{\Delta}$. In particular, we have the following

\begin{proposition}\label{wit}
Let $Z:R_{\Delta}\to \CO$ be 
an operator commuting with the action of differential operators
with constant coefficients. Then we have, for $f\in G_{\Delta}$, 
$$
Z(f)(0)=Tr_{S_{\Delta}}\left(Res_{\Delta}m(f)C Z Res_{\Delta}\right)
$$
where $(CZ)(\phi)(z)=Z(\phi)(-z)$.
\end{proposition}

Choosing a diagonal basis of $O\CB(\Delta)$, we can express the
preceding formula as a residue formula in several variables:

$$
Z(f)(0)=\sum_{o\sigma\in OB}
Res^{o\sigma}(f(z)Z_{\sigma}(-z))
$$
with $Z_{\sigma}(z)=Z(\phi_{\sigma})(z)$.

For later use, we prove a vanishing property 
of the linear form $Res^{o\sigma}$. Let $o\sigma$ be an ordered
basis. We write $o\sigma=(\alpha_1,\alpha_2,\ldots,\alpha_r)$ and 
$z=(z_1,z_2,\ldots,z_r)$. Set $o\sigma'=(\alpha_2,\ldots,\alpha_r)$
and $z'=(z_2,\ldots,z_r)$; then $z=(z_1,z')$. Let 
$\psi(z')$ in $\CO_{\Delta'}$ be a meromorphic function with denominator
a product of linear forms $\alpha(z')$ where $\alpha\in\Delta$ is not
a multiple of $\alpha_1$.

\begin{lemma}\label{b}
For any $f\in G_{\Delta}$,
and any $\psi\in \CO_{\Delta'}$,
$$
Res^{o\sigma}\left(\frac{1}{z_1}f(z_1,z')\psi(z')\right)=0.
$$
\end{lemma}

\begin{proof}
We have 
$$
Res^{o\sigma}\left(\frac{1}{z_1}f(z_1,z')\psi(z')\right)=
Res_{z_1=0}\left(\frac{1}{z_1}Res^{o\sigma'}(f(z_1,z')\psi(z'))\right).
$$
In computing $Res^{o\sigma'}\left(f(z_1,z')\psi(z')\right)$,
the variable $z_1$ is fixed to a non-zero value.
The result $Res^{o\sigma'}\left(f(z_1,z')\psi(z')\right)$
is a meromorphic function of $z_1$. It is thus sufficient to prove that 
$Res^{o\sigma'}\left(f(z_1,z')\psi(z')\right)$ belongs to the space 
$G=\sum_{k\leq -1}\C z_1^k$. 

We check this for $f=\phi_{\kappa}$ where  
$$
\phi_{\kappa}(z)=\frac{1}{\prod_{\alpha\in \kappa}\langle \alpha,z\rangle}
$$
and $\kappa$ is a generating sequence. Let 
$$
\kappa_1:=\{\alpha\in \kappa, \langle\alpha,(z_1,0)\rangle\neq 0\}
$$
and 
$$
\kappa'=\{\alpha\in \kappa, \langle\alpha,(z_1,0)\rangle=0\}.
$$
As $\kappa$ is generating, the set $\kappa_1$ is not empty.
We fix $z_1\neq 0$. We have
$$
\phi_{\kappa}(z_1,z')\psi(z')=
\phi_{\kappa_1}(z_1,z')\phi_{\kappa'}(z')\psi(z')
$$
and $\phi_{\kappa'}\in\CO_{\Delta'}$.
For $\alpha \in \kappa_1$, we set
$\langle \alpha,(z_1,z')\rangle =c_\alpha z_1+\langle \beta,z'\rangle$, 
with $c_\alpha\neq 0$. We consider the Taylor expansion at $z'=0$
of the holomorphic function of $z'$:
$$
\frac{1}{\langle \alpha, (z_1,z')\rangle }=
\frac{1}{c_\alpha z_1+\langle \beta,z'\rangle }=\frac{1}
{c_\alpha z_1(1+\frac{\langle \beta,z'\rangle }{c_\alpha z_1})}.
$$
This is of the form 
$$
\sum_{k=1}^{\infty} z_1^{-k}P_{k-1}(z')
$$
where $P_{k-1}(z')$ is homogeneous of degree $k-1$ in $z'$.
Let $n=|\kappa_1|$, then $n\geq 1$. We see that the function
$$z'\mapsto \phi_{\kappa_1}(z_1,z')=
\frac{1}{\prod_{\alpha\in \kappa_1}\langle \alpha, (z_1,z')\rangle }
$$
has a Taylor expansion of the form 
$$
\sum_{k\geq n} z_1^{-k}Q_{k-1}(z')
$$
where $Q_{k-1}(z')$ is homogeneous of degree $k-1$ in $z'$. Thus
$$
Res^{o\sigma'}
\left(\phi_{\kappa_1}(z_1,z')\phi_{\kappa'}(z')\psi(z')\right)
= \sum_{k\geq n} z_1^{-k} Res^{o\sigma'}
\left( Q_{k-1}(z')\phi_{\kappa'}(z')\psi(z')\right).
$$
Via the Taylor series at $z'=0$, the function 
$\phi_{\kappa'}(z')\psi(z')$ can be expressed as an infinite sum 
of homogeneous elements with finitely many negative degrees. 
As the iterated residue $Res^{o\sigma'}$ vanishes on elements of degree 
not equal to $-(r-1)$, and $Q_{k-1}(z')$ is homogeneous of degree $k-1$,
we see that the sum is finite and that
$Res^{o\sigma'}
\left(\phi_{\kappa_1}(z_1,z')\phi_{\kappa'}(z')\psi(z')\right)$
is in the space $G$ as claimed.
\end{proof}

\section{Eisenstein series}

Results of the second section will be used for
a complex vector space which is the complexification of a real 
vector space. Thus we slightly change the notation in this section.

Let $V$ be a {\bf real} vector space of dimension $r$ equipped with 
a lattice $N$. The complex vector space $V_\C$ is the space to
which we will apply the results of Section 2.
 
We consider the dual lattice $M=N^*$ to $N$.
We consider the  compact torus $T=iV/(2i\pi N)$
and its complexification $T_\C=V_\C/(2i\pi N)$.
The projection map $V_\C\to T_\C$ is denoted by the exponential notation
$v\to e^v$. If $\{e^1, e^2,\ldots,e^r\}$ is a $\Z$-basis of $N$, we
write an element of $V_\C$ as $z= z_1 e^1+z_2 e^2+ \cdots +z_r e^r$ 
with $z_j\in \C$. We can identify 
$T_\C$ with $\C^*\times \C^*\times \cdots \times \C^*$
by $z\mapsto (e^{z_1}, e^{z_2},\ldots,e^{z_r})$.

If $m\in M$, we denote by $e^m$ the character of $T$ defined by 
$\langle e^m,e^v\rangle = e^{\langle m,v\rangle}$.
We extend $e^m$ to a holomorphic character of the complex torus $T_{\C}$.
The ring of holomorphic functions on $T_\C$ generated by the functions
$e^m$ is denoted by $R(T)$. A quotient of two elements of 
$R(T)$ is called a {\bf rational function} on the complex torus $T_\C$. 
Via the exponential map $V_\C\to T_\C$, a function on $T_\C$
will be sometimes identified with a function 
on $V_\C$, invariant under translation by the lattice $2i\pi N$. 
If $\{e^1, e^2,\ldots,e^r\}$ is a $\Z$-basis of $N$, a rational function 
on $T_\C$ written in exponential coordinates
is a rational function of $e^{z_1}, e^{z_2},\ldots,e^{z_r}$.
We say briefly that it is a rational function of $e^z$.

Let us consider a finite set $\Delta$ of non-trivial characters of 
$T$. We identify $\Delta$ with a subset of $M$; for $\alpha\in\Delta$,
we denote by  $e^{\alpha}$ the corresponding character of $T_\C$.
\begin{definition}
We denote by $R(T)_{\Delta}$ the subring of rational functions on $T$
generated by $R(T)$ and the inverses of the functions 
$1-e^{-\alpha}$ with $\alpha\in \Delta$.
\end{definition}

Observe that $R_{\Delta}$ is left unchanged when each element of
$\Delta$ is replaced by a non-zero scalar multiple, but that
$R(T)_{\Delta}$ strictly increases when (say) each
$\alpha\in\Delta$ is replaced by $2\alpha$. We assume from now on 
that all elements of $\Delta$ are indivisible in the lattice $M$.

Via the exponential map, we consider 
elements of $R(T)_{\Delta}$
as periodic meromorphic functions on $V_{\C}$.
On $V_\C$, the function 
$$
\frac{\langle \alpha,z\rangle }{1-e^{-\langle \alpha,z\rangle }}
$$
is defined at $z=0$, so is an element of $\CO$.
Writing 
$$
\frac{1}{1-e^{-\langle \alpha,z\rangle }}=
\frac{1}{\langle \alpha,z\rangle}
\frac{\langle \alpha,z\rangle }{1-e^{-\langle \alpha,z\rangle }},
$$
we see that $R(T)_{\Delta}$ is contained in $\CO_{\Delta}$.
We see furthermore from the formula 
$$
\frac{d}{dz}\frac{1}{1-e^{-z}}=
\frac{1}{(1-e^z)(1-e^{-z})}=
\frac{-e^{-z}}{(1-e^{-z})^2 }
$$
that $R(T)_\Delta\subset\CO_{\Delta}$ is stable under differentiation.

Our aim is to find a natural map from 
$R_{\Delta}$ to $R(T)_{\Delta}$ commuting with the action of 
differential operators with constant coefficients.
In particular, we want to force a rational function of $z\in V_\C$
to become periodic, so that it is natural to define Eisenstein series 
$$
E(f)(z)=\sum_{n\in N} f(z+2i\pi n).
$$
We need to be more careful, as the sum is usually not convergent
for an arbitrary $f\in R_{\Delta}$.
We will introduce an oscillating factor $e^{\langle t,2i\pi n\rangle }$
with $t\in V^*$ in front of each term of this infinite sum.

Let 
$$
U_{\Delta}=\{z\in V_\C,\langle \alpha,z+2i\pi n\rangle \neq 0
~{\rm for~all~}n\in N~{\rm and~for~all}~ \alpha\in \Delta\}.
$$
Then $R(T)_{\Delta}$ consists of periodic holomorphic functions on
$U_{\Delta}$.

Let $f\in R_{\Delta}$, then $f(z+2i\pi n)$ is defined for each 
$n\in N$ if $z\in U_\Delta$. For $z\in U_{\Delta}$, we consider 
the function on $V^*$ defined by  
$$
t\mapsto 
\sum_{n\in N}e^{\langle t, z+2i\pi n\rangle }f(z+2i\pi n).
$$
If $n\mapsto f(z+2i\pi n)$
is sufficiently decreasing at infinity,
the series will be absolutely convergent and sum up to
a continuous function of $t$ with value at $t=0$ equal to 
$$
\sum_{n\in N}f(z+2i\pi n).
$$
In any case, it is easy to see that this series of functions 
of $t$ converges to a generalized function of $t$:

\begin{proposition}
For each $f\in R_{\Delta}$ and $z\in U_{\Delta}$,
the function on $V^*$ defined by  
$$
t\mapsto 
\sum_{n\in N}e^{\langle t, z+2i\pi n\rangle }f(z+2i\pi n)
$$
is well defined as a generalized function of $t$, which depends
holomorphically on $z$ for $z$ in the open set $U_\Delta$.
\end{proposition}
\begin{proof}
Indeed, if $s(t)$ is a smooth function on $V^*$ with compact support, 
consider the series
$$
\sum_{n\in N}f(z+2i\pi n)\int_{V^*} e^{\langle t,z+2i\pi n\rangle}s(t)dt
=\sum_{n\in N}c(z,n)f(z+2i\pi n).
$$
The coefficient
$$
c(z,n)=\int_{V^*}e^{2i\pi\langle t,n\rangle}e^{\langle t,z\rangle}s(t)dt
$$
is rapidly decreasing in $n$, as 
the function $t\mapsto e^{\langle t,z\rangle}s(t)$ is smooth 
and compactly supported. Thus $c(z,n) f(z+2i\pi n)$
is also a rapidly decreasing function of $n$. Furthermore 
$c(z,n) f(z+2i\pi n)$ depends holomorphically on $z\in U_{\Delta}$. So
the result of the summation 
$$
\sum_{n\in N}c(z,n) f(z+2i\pi n)
$$
exists and is a holomorphic function of $z$.
\end{proof}

We write 
$$
E(f)(t,z)=\sum_{n\in N}e^{\langle t, z+2i\pi n\rangle }f(z+2i\pi n)
$$
for this generalized function of $t$ depending holomorphically 
on $z$. We will analyze this function of $(t,z)$, $t\in V^*$, 
$z\in U_{\Delta}$.

We summarize first some of the obvious properties of $E(f)(t,z)$. 
\begin{proposition}\label{dif}
The following equations are satisfied:

1) For every $P\in S(V^*)$ and $f\in R_{\Delta}$,
$$
E(Pf)(t,z)=P(\partial_t) E(f)(t,z).
$$

2) For every $v\in V$ and $f\in R_{\Delta}$,
$$
E(\partial(v)f)(t, z)
=\partial_z(v) E(f)(t,z)-\langle t,v\rangle E(f)(t,z).
$$

3) For every  $m\in M$ and $z\in U_{\Delta}$,
$$
E(f)(t+m,z)=e^{\langle m,z \rangle }E(f)(t,z).
$$

\end{proposition}
As $R_\Delta$ is generated by $S_\Delta$ under the action of 
$S(V)$ and $S(V^*)$, we see that the operator 
$E$ is completely determined by the functions $E(\phi_{\sigma})(t,z)$
($\sigma\in\CB(\Delta)$). 

A {\bf wall} of $\Delta$ is a hyperplane of 
$V^*$ generated by $r-1$ linearly independent vectors of $\Delta$.
We consider the system of affine hyperplanes generated 
by the walls of $\Delta$ together with their translates 
by $M$ (the dual lattice of  $N$). We denote by $V^*_{\Delta,areg}$ 
the complement of the union of these affine hyperplanes.
A connected component of $V^*_{\Delta,areg}$ will be called an 
{\bf alcove} and will be denoted by ${\mathfrak a}$.

\begin{proposition}\label{pol}
The function $E(f)(t,z)$ is smooth when $t$ varies on 
$V^*_{\Delta,areg}$, and $z\in U_{\Delta}$. More precisely, let
${\mathfrak a}$ be an alcove. Assume that $f$ is homogeneous 
of degree $d$. Then, on the open set ${\mathfrak a}\times U_\Delta$, the
function $E(f)(t,z)$ is a polynomial in $t$ of degree at most $-d-r$,
with coefficients in $R(T)_{\Delta}$. 
\end{proposition}
\begin{proof}
Consider first the one variable case. The set 
$V^*_{\Delta,areg}$ is $\R-\Z$.
Let $[t]$ be the integral part of $t$. Fix $z\in\C - 2i\pi\Z$.
Consider the locally  constant function of $t\in\R-\Z$ defined by 
$$
t\mapsto \frac{e^{[t]z}}{1-e^{-z}}.
$$
We extend this function as a locally $L^1$-function on $\R$
(defined except on the set $\Z$  of measure $0$).
\begin{lemma}
We have the equality of generalized functions of $t$:
$$
\sum_{n\in\Z}\frac{e^{t(z+2i\pi n)}}{z+2i\pi n}
=\frac{e^{[t]z}}{1-e^{-z}}.
$$
\end{lemma}
\begin{proof}
We compute the derivative in $t$ of the left hand side.
It is equal to
$$
\sum_{n\in \Z}e^{t(z+2i\pi n)}=e^{tz}\delta_\Z(t)
$$
where $\delta_\Z$ is the delta function of the set of integers.

We compute the derivative in $t$ of the right hand side.
This function of $t$ is constant on each interval $(n, n+1)$.
The jump at the integer $n$ is 
$$
\frac{e^{nz}}{1-e^{-z}}-\frac{e^{(n-1)z}}{1-e^{-z}}
=e^{nz}.
$$
It follows that the derivative in $t$ of the right hand side is also
equal to $e^{tz}\delta_\Z(t)$. Thus  
$$
\sum_{n\in\Z}\frac{e^{t(z+2i\pi n)}}{z+2i\pi n}=
c(z)+\frac{e^{[t]z}}{1-e^{-z}}
$$
where $c(z)$ is a constant.
We verify that $c(z)$ is equal to $0$ by using periodicity properties
in $t$. It is clear that
$$
e^{-tz} \sum_{n\in\Z}\frac{e^{t(z+2i\pi n)}}{z+2i\pi n}=
\sum_{n\in\Z}\frac{e^{2i\pi n t}}{z+2i\pi n}
$$
is a periodic function of $t$ as is 
$$
e^{-tz}\frac{e^{[t]z}}{1-e^{-z}}= 
\frac{e^{([t]-t)z}}{1-e^{-z}}.
$$
It follows that $e^{-tz}c(z)$
is also a periodic function of $t$. This implies $c(z)=0$.
\end{proof}

Consider now, for $k\in\Z$,
$$
E_k(t,z)=\sum_{n\in\Z}e^{t(z+2i\pi n)}(z+2i\pi n)^k.
$$
We just saw that
$$
E_{-1}(t,z)=\frac{e^{[t]z}}{1-e^{-z}}.
$$
To determine $E_k(t,z)$ for $k\leq -1$, we use the differential 
equation in $z$ 
$$
\partial_z E_k(t,z)=t E_k(t,z) + k E_{k-1}(t,z).
$$
Using decreasing induction over $k$, we see that $E_{k}(t,z)$ is a
$L^1$-function of $t$, equal to a polynomial function 
of $t$ of degree $-k-1$ on each interval $(n,n+1)$
and with coefficients rational functions of $e^z$.
For example, we obtain the value of the convergent series
$$
\sum_n \frac{e^{t(z+2i\pi n)}}{(z+2i\pi n)^2}=
(t-[t])\frac{e^{[t]z}}{1-e^{-z}}-\frac{e^{[t]z}}{(1-e^{-z})(1-e^{z})}.
$$
When $k\geq 0 $, we use the differential equation 
$$
\partial_t E_k(t,z)=E_{k+1}(t,z)
$$  
so that, as we have already used,
$$
E_0(t,z)=\sum_{n\in \Z}e^{t(z+2i\pi n)}=e^{tz}\delta_\Z(t).
$$
More generally, $E_k(t,z)=(\partial_t)^k(e^{tz}\delta_\Z(t))$ 
is supported on $\Z$, in particular is identically $0$ on $\R-\Z$.
\bigskip

We return to the proof of Proposition \ref{pol}. For a simple fraction
$\phi$, consider the function 
$$
t\mapsto E(\phi)(t, z).
$$
We first prove that it is a locally $L^1$-function, which is constant
when $t$ varies in an alcove.

Let $\sigma=\{ \alpha_1, \alpha_2,\ldots, \alpha_r\}$ be a basis of
$\Delta$. Let $t\in V^*$. If $t=\sum_j t_j\alpha_j$ is the decomposition 
of $t$ on the basis $\sigma$, set $[t]_{\sigma}=\sum_j [t_j]\alpha_j$. The
function $t\mapsto [t]_{\sigma}$ is constant when $t$ varies in an alcove.
Consider the sublattice 
$$
M_{\sigma}=\bigoplus_{\alpha\in \sigma} \Z \alpha \subseteq M.
$$ 
We say that $\sigma$ is a $\Z$-basis, if $M_{\sigma}=M$. In general,
the quotient $M/M_{\sigma}$ is a finite set; let ${\cal R}$ be a set
of representatives of this quotient. We can choose ${\cal R}$ in the
following standard way. We consider the box 
$$
Q_{\sigma}=\bigoplus_{\alpha\in \sigma} [0,1[\alpha
= \{u\in V^*, [u]_{\sigma}=0\}.
$$
Then we can take
$$
{\cal R}=Q_{\sigma}\cap M=\{u\in M,[u]_{\sigma}=0\}.
$$
Define 
$$
{\cal R}(t,\sigma)=(t-Q_{\sigma})\cap M=\{u\in M, [t-u]_{\sigma}=0\}.
$$
The set ${\cal R}(t,\sigma)$ is also a set of representatives of
$M/M_{\sigma}$. If $\sigma$ is a $\Z$-basis of $M$, 
this set is reduced to the single element $[t]_{\sigma}$.  
Remark that the set ${\cal R}(t,\sigma)$
is constant when $t$ varies in an alcove ${\mathfrak a}$.
We denote it by ${\cal R}({\mathfrak a},\sigma)$. 
\begin{definition}\label{eis}
If $\mathfrak a$ is an alcove, and $\sigma$ a basis of $\Delta$, we set
$$
F_{\sigma}^{{\mathfrak a}}=|M/M_{\sigma}|^{-1}
\frac{\sum_{m\in {\cal R}({{\mathfrak a},\sigma)}}e^m}
{\prod_{\alpha\in \sigma}(1-e^{-\alpha})}.
$$
\end{definition}
Thus an alcove ${\mathfrak a}$ together with a basis
$\sigma\in\CB(\Delta)$ produces a particular element
$F_{\sigma}^{{\mathfrak a}}$ of $R(T)_{\Delta}$.

Consider on the set $V^*_{\Delta,areg}$
the locally constant function of $t$ defined by 
$F_{\sigma}(t,z)=F_{\sigma}^{{\mathfrak a}}(z)$ when $t$ is in the
alcove ${\mathfrak a}$.
This defines a locally $L^1$-function of $t$, still denoted 
by $F_{\sigma}(t,z)$, defined except on the set 
$V^*-V^*_{\Delta,areg}$ 
of measure $0$. This locally $L^1$-function of $t$ defines a
generalized function of $t$ which depends holomorphically on $z$.

\begin{lemma}
We have the equality of generalized functions of $t\in V^*$:
$$
E(\phi_{\sigma})(t,z)=F_{\sigma}(t,z).
$$
\end{lemma}
\begin{proof}
If $\sigma$ is a $\Z$-basis of $M$, this follows from 
the formula in dimension $1$. In general, we consider 
$M_{\sigma}\subseteq M$ and the dual lattice $N_{\sigma}=M_{\sigma}^*$.
Then $N\subseteq N_{\sigma}$. We set
$$
E_{\sigma}(\phi_{\sigma})(t,z):=\sum_{\ell\in N_{\sigma}}
e^{\langle t,z+2i\pi\ell\rangle}\phi_{\sigma}(z+2i\pi \ell).
$$
For any set of representatives ${\cal R}$ 
of $M/M_{\sigma}$, we have:
$\sum_{u\in{\cal R}}e^{-\langle u,2i\pi \ell\rangle}=0$
if $\ell\in N_{\sigma}$ is not in $N$,
while this sum equal $|M/M_{\sigma}|$ if $n\in N$. Thus,
$$\displaylines{
E(\phi_{\sigma})(t, z)=
\sum_{n\in N}\phi_{\sigma}(z+2i\pi n)e^{\langle t,z+2i\pi n\rangle }
\hfill\cr
= \sum_{\ell\in N_{\sigma}}
\phi_{\sigma}(z+2i\pi\ell) e^{\langle t,z+2i\pi \ell\rangle }
(|M/M_{\sigma}|^{-1}\sum_{u\in{\cal R}}e^{-\langle u,2i\pi\ell\rangle})
\cr
=|M/M_{\sigma}|^{-1}\sum_{u\in{\cal R}}\sum_{\ell\in N_{\sigma}}
\phi_{\sigma}(z+2i\pi\ell)
e^{\langle t-u,z+2i\pi \ell\rangle }e^{\langle u,z\rangle}
\cr\hfill
=|M/M_{\sigma}|^{-1}\sum_{u\in{\cal R}}e^{\langle u,z\rangle} 
E_{\sigma}(\phi_{\sigma})(t-u,z).
\cr}$$
This holds as an equality of generalized functions of $t$.
Further, we have by the one-dimensional case:
$$
E_{\sigma}(\phi_{\sigma})(t,z)=\frac{e^{\langle [t]_{\sigma},z\rangle }}
{\prod_{\alpha\in \sigma}(1-e^{-\langle \alpha,z\rangle })}.
$$
It follows that $E(\phi_{\sigma})(t,z)$ is a locally $L^1$-function of
$t$, as is $E_{\sigma}(\phi_{\sigma})$. It remains to determine the
value of this function when $t$ is in an alcove. For $m\in M_{\sigma}$,
we have
$$
E_{\sigma}(\phi_{\sigma})(t+m,z)=
e^{\langle m, z\rangle }E_{\sigma}(\phi_{\sigma})(t,z),
$$
so that the sum 
$\sum_{u\in{\cal R}}e^{\langle u,z\rangle} 
E_{\sigma}(\phi_{\sigma})(t-u,z)$
is independent of the choice of the system 
of representatives ${\cal R}$ of $M/M_{\sigma}$. We choose
${\cal R}={\cal R}(t,\sigma)$. Then 
$$
E(\phi_{\sigma})(t,z)=|M/M_{\sigma}|^{-1}
\frac{\sum_{u\in {\cal R}(t,\sigma)} e^{\langle u,z\rangle }}
{\prod_{\alpha\in \sigma}(1-e^{-\langle \alpha,z\rangle })}
$$
because $[t-u]_{\sigma}=0$ for all $u\in {\cal R}(t,\sigma)$.
\end{proof}

Every function $f\in R_{\Delta}$, homogeneous of degree $d$, is
obtained from an element of $S_\Delta$ by the action 
of a differential operator with polynomial coefficients. 
This operator is of degree $d+r$, if multiplication 
by $z_j$ is given degree $1$, while derivation 
$\frac{\partial}{\partial z_j}$ is given degree $-1$.
Using Proposition \ref{dif}, we see that Proposition \ref{pol}
follows from the fact that the function
$t\mapsto E(\phi_{\sigma})(t,z)$ is constant on each alcove.
\end{proof}

\bigskip

From Proposition \ref{pol}, we see that 
there exist functions $\phi_{(k)}^{\mathfrak a}(z)\in R(T)_{\Delta}$ 
such that we have the equality  for $t$ in the alcove 
${\mathfrak a}$:
$$
E(f)(t,z)= \sum_{n\in N}e^{\langle t, z+2i\pi n\rangle }f(z+2i\pi n)
=\sum_{(k)}t^{(k)}\phi^{\mathfrak a}_{(k)}(z)$$
where the sum is over a finite number of multi-indices $(k)$. This
defines an operator 
$$
E^t:R_\Delta\to R(T)_{\Delta}, f\mapsto E(f)(t,z)
$$
obtained by fixing the regular value $t$.

The operator $E^t$ satisfies the following relation, which is just
the relation $(2)$ in Proposition \ref{dif}: For $v\in V$ and  
$f\in R_{\Delta}$,
$$ 
E^t(\partial(v)f)(z)=
\partial_z(v) E^t(f)(z)-\langle t,v\rangle E^t(f)(z).
$$

Let $B$ be a basis of $\CB(\Delta)$.
Let ($\phi_{\sigma}$, $\sigma\in B$) be the corresponding basis of 
$S_{\Delta}$ and ($\phi^{\sigma}$, $\sigma\in B$) the dual basis of
$S_{\Delta}^*$. For $\sigma\in B$, and an alcove ${\mathfrak a}$,
consider the element 
$F_{\sigma}^{{\mathfrak a}}$ of $R(T)_{\Delta}\subset \CO_{\Delta}$
associated to $\sigma,{\mathfrak a}$. We obtain a kernel formula for
the operator $E^t$:

\begin{theorem}\label{ke}
Let $f\in G_{\Delta}$. For $y\in U_{\Delta}$ and $t\in{\mathfrak a}$,
we have:
$$\displaylines{
E^t(f)(y)=Tr_{S_{\Delta}}
\left( Res_{\Delta} m(e^{\langle t,\cdot\rangle}f) C
{\cal T}(y) E^t Res_{\Delta}\right)
\hfill\cr\hfill
=\sum_{\sigma\in B} 
\langle \phi^\sigma,Res_{\Delta}\left(e^{\langle t, z\rangle}f(z)
F_{\sigma}^{\mathfrak a}(y-z)\right)\rangle
\cr}$$
where $F_{\sigma}^{\mathfrak a}$ is given by Definition \ref{eis}. If
moreover $B$ is the underlying basis of a diagonal basis $OB$, then 
$$
E^t(f)(y)=\sum_{o\sigma\in OB} Res^{o\sigma}
\left(e^{\langle t,z\rangle} f(z) F_{\sigma}^{\mathfrak a}(y-z)\right).
$$
\end{theorem}
\begin{proof}
By a method entirely similar to the proof of Theorem \ref{ker}, we see
that the operator
$$
A^t(f)(y)=\sum_{\sigma\in B} 
\langle \phi^{\sigma}, Res_{\Delta}\left(e^{\langle t,z\rangle}f(z)
F_{\sigma}^{\mathfrak a}(y-z)\right)\rangle 
$$
satisfies the relation
$$
A^t(\partial(v)f)( z)
=\partial_z(v) A^t(f)(z)-\langle  t,v\rangle A^t(f)(z)
$$
for $v\in V$, $f\in R_{\Delta}$. Thus to prove that 
$E^t=A^t$ on $G_{\Delta}$, it is sufficient to prove that 
they coincide for $f=\phi_{\tau}$. In this case, we obtain 
$$
A^t(\phi_{\tau})(y) = \sum_{\sigma\in B} \langle \phi^{\sigma}, 
\phi_{\tau}(z)\rangle F_{\sigma}^{\mathfrak a}(y) =
F_{\tau}^{\mathfrak a}(y)=E^t(\phi_{\tau})(y).
$$ 
\end{proof}

In view of the kernel formula for the Eisenstein series 
$E^t$, it is natural to introduce the following definition.

\begin{definition}
The {\bf constant term} of the Eisenstein series $E^t$ is the linear
form $f\to CT(f)(t)$ defined for $f\in R_{\Delta}$ 
and $t$ in the alcove ${\mathfrak a}$ by 
$$
CT(f)(t)= 
Tr_{S_\Delta}\left(Res_{\Delta}m(e^{\langle t,\cdot\rangle}f)C E^t 
Res_{\Delta}\right).
$$
More explicitly, if $OB$ is a diagonal basis of $\CB(\Delta)$, then
$$
CT(f)(t)=\sum_{o\sigma\in OB} 
Res^{o\sigma}\left(e^{\langle t,z\rangle}f(z)
F_{\sigma}^{{\mathfrak a}}(-z)\right).
$$
\end{definition}

\section{Partial Eisenstein series}

Let $N_{reg}=N\cap V_{reg}$ be the set of regular elements of $N$. The
aim of this section is to prove that the function
$$
E_{N_{reg}}(f)(t,z)=
\sum_{n\in N_{reg}} e^{\langle t,z+2i\pi n\rangle} f(z+2i\pi n)
$$
is analytic in $(t,z)$ when $t$ is in an alcove and $z\in V_{\C}$ is
close to $0$. In the next section we will prove Szenes residue formula
for
$$
E_{N_{reg}}(f)(t,0)=
\sum_{n\in N_{reg}}e^{\langle t,2i\pi n\rangle} f(2i\pi n).
$$

Let $\Gamma$ be a subset of $N$. We can define, for $f\in R_{\Delta}$, 
the generalized function of $t$ 
$$
E_{\Gamma}(f)(t,z)=\sum_{n\in \Gamma}
e^{\langle t, z+2i\pi n\rangle} f(z+2i\pi n).
$$
Introduce the set 
$$
U_{\Delta,\Gamma}=\{z\in V_{\C},\langle\alpha, z+2i\pi n\rangle\neq 0
~{\rm for~all}~\alpha\in \Delta~{\rm and}~n\in \Gamma\}.
$$
The generalized function $E_{\Gamma}(f)(t,z)$ depends holomorphically
on $z$, when $z\in U_{\Delta,\Gamma}$.
 
Let $W$ be a rational subspace of $V$. Then $N\cap W$ is a lattice in $W$.
Consider, for $f\in R_{\Delta}$,
$$
E_{N\cap W}(f)(t,z)=
\sum_{n\in N\cap W}e^{\langle t, z+2i\pi n\rangle} f(z+2i\pi n).
$$
We analyze the singularities in $(t,z)$ of $E_{N\cap W}(f)(t,z)$.
If $W$ is zero, then $E_{\{0\}}(f)(t,z)=e^{\langle t,z\rangle}f(z)$ is
analytic in $(t,z)$ when $z$ is regular in $V_{\C}$. Assume that $W$
is non-zero and consider the subspace $W^{\perp}$ of $V^*$.
Remark that, if $u\in M+W^{\perp}$, we have the relation 
$$
E_{N\cap W}(f)(t+u,z)=e^{\langle u,z\rangle }E_{N\cap W}(f)(t,z).
$$
It is clear that the singular set of $E_{N\cap W}(f)(t,z)$
is stable by translation by $M+W^{\perp}$. Define a 
$(W,\Delta)$-{\bf wall} in $V^*$ as an hyperplane generated by
$W^{\perp}$ together with $\dim W-1$ vectors of $\Delta$.
We introduce the set $\CH^*_{W,\Delta,M}$ consisting of the union of all
$(W,\Delta)$-walls and of their translates by elements of $M$.
We define $V^*_{W,\Delta,areg}$ as the complement of
$\CH^*_{W,\Delta,M}$ in $V^*$. This set $V^*_{W,\Delta,areg}$ is
invariant by translation by $M+W^{\perp}$. 

\begin{lemma}\label{par}
For $f\in R_{\Delta}$, the function $E_{N\cap W}(f)(t,z)$ is
analytic in $(t,z)$ when $t$ varies on $V^*_{W,\Delta,areg}$ and
$z\in U_{\Delta,N\cap W}$. Furthermore, if $t\in V^*_{W,\Delta,areg}$
and $z$ is near $0$, the function $z\mapsto E_{N\cap W}(f)(t,z)$
defines an element of $\CO_{\Delta}$.
\end{lemma}
\begin{proof}
Let $\sigma$ be a basis of $\Delta$. Although we are not able to give
a nice formula for the function 
$E_{N\cap W}(\phi_{\sigma})(t,z)$, 
we can still obtain an inductive expression that suffices to give
some informations on it. Consider the set $V^*_{W,\sigma,areg}$, that
is, the complement of the union of $(W,\sigma)$-walls together with
their translates by $M$. Let $U_{\sigma,N\cap W}$ be the set of all
$z\in V_{\C}$ such that $\langle \alpha, z+2i\pi n\rangle \neq 0$
for all $\alpha\in \sigma$ and $n\in N\cap W$.
The intersection of this set with a small neighborhood of $0$ is
contained in the complement of the union of the complex hyperplanes 
$\{z\in V_{\C},\langle \alpha,z\rangle=0\}$,
for $\alpha\in \sigma$.

\begin{lemma}\label{ana}
The function $E_{N\cap W}(\phi_{\sigma})(t,z)$ is
analytic in $t\in V^*_{W, \sigma,areg}$ and $z\in U_{\sigma,N\cap W}$.
Furthermore, when $t\in  V^{*}_{W, \sigma, areg}$, the function 
$$
z\mapsto
\left(\prod_{\alpha\in \sigma}\langle \alpha,z\rangle \right)
E_{N\cap W}(\phi_{\sigma})(t,z)
$$ 
is holomorphic at $z=0$.
\end{lemma}
We prove this by induction on the codimension of $W$.
If $W=V$, this follows from the explicit formula for
$E(\phi_{\sigma})(t,z)$. Let $\alpha$ be an indivisible element of $M$
such that $W$ is contained in the real hyperplane 
$$
H_{\alpha}=\{y\in V,\langle\alpha,y\rangle=0\}.
$$
We assume first that $\alpha$ is an element of $\sigma$.
We number it the first vector $\alpha_1$ of the basis $\sigma$.
We set $\sigma'=(\alpha_2,\ldots,\alpha_r)$, $z'=(z_2,\ldots,z_r)$,
etc; then $z=(z_1,z')$. Our subspace $W$ is contained in
$V'=V\cap\{z_1=0\}$. Thus, we have
$$
E_{N\cap W}(\phi_{\sigma})(t,z)=
\sum_{n\in N\cap W}e^{\langle t, z+2i\pi n\rangle} \phi_{\sigma}(z+2i\pi n)
=\frac{e^{t_1 z_1}}{z_1}E_{N'\cap W}(\phi_{\sigma'})(t',z').
$$
By induction, $E_{N'\cap W}(\phi_{\sigma'})(t',z')$
is analytic in $(t',z')$ for $z'\in U_{\sigma',N'}$,  
except if there exist $m'\in M'$ such that 
$t'+m'$ is in a hyperplane generated by $W^{\perp'}$
(the orthogonal of $W$ in $V'$) and some vectors of $\sigma'$. 
As $W^{\perp}=W^{\perp'}\oplus \R \alpha_1$,
we see that the singular set of $E_{N\cap W}(\phi_{\sigma})(t,z)$
is contained in $\CH^*_{W,\sigma,M}$.
Furthermore, the function 
$$
z_1z_2\cdots z_r E_{N\cap W}(\phi_{\sigma})(t,z)= e^{t_1 z_1}
z_2\cdots z_r E_{N'\cap W}(\phi_{\sigma'})(t',z')
$$
is holomorphic in $z$ near $z=0$.

Assume now that $\alpha$ is not an element of $\sigma$.
We add it to the system $\Delta$ if $\alpha$ is not an element 
of $\Delta$.
Writing $\alpha=\sum_j c_j \alpha_j$, we obtain one of the 
Orlik-Solomon relations of the system $\Delta\cup \{\alpha\}$ 
$$
\phi_{\sigma}=\sum_j c_j \phi_{\sigma^j}
$$
where ${\sigma}^j=\sigma\cup\{\alpha\}-\{\alpha_j\}$. A
$(W,\sigma^j)$-wall is a hyperplane of $V^*$
generated by $W^{\perp}$ and $\dim W -1$ vectors of $\sigma^j$; then
these vectors are distinct from $\alpha$, because 
$\alpha\in W^{\perp}$. Thus, all $W$-walls for the basis $\sigma^j$ are
also $W$-walls for the basis $\sigma$. By our first calculation, it
follows that $E_{N\cap W}(\phi_{\sigma^j})(t,z)$ is analytic when
$t$ is not on a translate of a $(W,\sigma)$-wall. Moreover, we have 
$$
E_{N\cap W}(\phi_{\sigma})(t,z)=
\sum_j c_j E_{N\cap W}(\phi_{\sigma^j})(t,z),
$$
so that the function 
$$
z\mapsto\langle \alpha,z\rangle 
\left(\prod_{j=1}^r \langle \alpha_j,z\rangle\right)
E_{N\cap W}(\phi_{\sigma})(t,z)
$$ 
is holomorphic in $z$ in a neighborhood of $0$.

By the induction hypothesis applied to $W\subseteq V'=\{\alpha=0\}$, the
function $z\mapsto E_{N\cap W}(\phi_{\sigma})(t,z)$ is holomorphic on
a non-empty open subset of $V'_{\C}$. So this function, considered as a
function of $z\in V_{\C}$, has no pole along $\alpha=0$.
This proves Lemma \ref{ana}, and hence Lemma \ref{par} when $f$ is a
simple fraction. The operator $E_{N\cap W}$ satisfies also the
commutation relation of Proposition \ref{dif}. Thus, using
differential operators with polynomial coefficients, we obtain the
statement of Lemma \ref{par} when $f$ is any element in $R_{\Delta}$. 
\end{proof}

Let $I$ be a subset of $\Delta$ and let 
$W_I=\cap_{i\in I}H_{\alpha_i}$. This is a rational subspace of $V$,
and the $(W_I,\Delta)$-walls are some of the walls of 
$\Delta$. Then it follows from Lemma 
\ref{par} that $E_{N\cap W_I}(f)(t,z)$
is a fortiori analytic when $t\in V^*_{areg}$
and $z\in U_{\Delta}$. 

\begin{definition}
A subset $\Gamma$ of $N$ is {\bf admissible}, if 
the characteristic function of $\Gamma$ is a linear combination 
of characteristic functions of sets 
$N\cap W_I$, where $I$ ranges over subsets of $\Delta$.
\end{definition}

Then we have by Lemma \ref{par}:

\begin{lemma}\label{admissible}
If $\Gamma$ is an admissible subset of $N$, the function 
$(t,z)\mapsto E_\Gamma(f)(t,z)$ is analytic when 
$t\in V^*_{\Delta,areg}$ and $z\in U_{\Delta,\Gamma}$. Furthermore,
when $z$ is near $0$ and $t\in V^*_{\Delta, areg}$, the function
$z\mapsto E_\Gamma(f)(t,z)$ defines an element of $\CO_{\Delta}$. 
\end{lemma}

If $\Gamma$ in an admissible subset of $N$, we can take 
the value at $t$ of the generalized function
$$
E_{\Gamma}(f)(t,z)=\sum_{n\in\Gamma}
e^{\langle t, z+2i\pi n\rangle } f(z+2i\pi n)$$
provided that $t$ is in an alcove ${\mathfrak a}$.
Thus, for $t\in {\mathfrak a}$, we can define the operator 
$E^t_\Gamma:R_{\Delta}\to \CO_{\Delta}, f\mapsto E_{\Gamma}(f)(t,z)$.
Now the argument of Theorem \ref{ke} proves 

\begin{proposition}
For $f\in G_{\Delta}$, $t\in V^*_{\Delta,areg}$ and 
$y\in U_{\Delta,\Gamma}$, we have
$$
E^t_\Gamma(f)(y)=
Tr_{S_\Delta}\left( Res_{\Delta}m(e^{\langle t,\cdot\rangle }f)
C{\cal T}(y)E^t_{\Gamma} Res_{\Delta}\right).
$$
More explicitly, if we choose a diagonal basis $OB$ then
$$
E^t_{\Gamma}(f)(y)=\sum_{o\sigma\in OB} Res^{o\sigma}
\left(f(z)e^{\langle t,z\rangle}F^t_{\Gamma,\sigma}(y-z)\right)
$$
where $F^t_{\Gamma,\sigma}(z)=E_{\Gamma}(\phi_{\sigma})(t,z)$.
\end{proposition}

\section{Witten series and Szenes formula}

For $f\in R_{\Delta}$, let us form the series 
$$
Z(f)(t,z)=\sum_{n\in N_{reg}} e^{\langle t,z+2i\pi n\rangle }
f(z+2i\pi n)
$$
where $N_{reg}$ is the set of regular elements of $N$.
Then $Z(f)(t,z)$ is defined as a generalized function 
of $t$. As $n$ varies in $N_{reg}$, this generalized function 
of $t$ depends holomorphically on $z$ when $z$ varies in 
a neighborhood of $0$. As $N_{reg}$ is an admissible subset of $N$, we
obtain from Lemma \ref{admissible}

\begin{proposition}
For any alcove ${\mathfrak a}$, $Z(f)(t,z)$ is an analytic
function of $(t,z)$ when $t\in {\mathfrak a}$ and $z$ is in a 
neighborhood of $0$.
\end{proposition}

We have
$$
Z(f)(t,0)=\sum_{n\in N_{reg}} e^{\langle t,2i\pi n\rangle}f(2i\pi n).
$$
This is well defined as a generalized function of $t$ when $t$ is in
an alcove. If $n\mapsto f(2i\pi n)$ is sufficiently decreasing, then
$Z(f)(t,0)$ is a continuous function of $t$; it generalizes the
Bernoulli polynomial 
$$
B_k(t)=\sum_{n\neq 0} \frac{e^{2i\pi n t}}{(2i\pi n)^k}
$$ 
where $0<t<1$.
 
We reformulate Szenes formula as an equality between 
$Z(f)(t,0)$ and the constant term of the Eisenstein series
$E(f)(t,z)$.

\begin{theorem}
For any $f\in R_\Delta$ and $t$ in an alcove ${\mathfrak a}$, we have 
$$
Z(f)(t,0)=CT(f)(t)= Tr_{S_\Delta}
\left(Res_{\Delta}m(e^{\langle t,\cdot\rangle} f)
C E^t Res_{\Delta}\right).
$$
In particular, $Z(f)(t,0)$ is a polynomial function of $t$ when $t$
varies in an alcove ${\mathfrak a}$.
\end{theorem}

As a consequence, if $OB$ is a diagonal basis, then we recover the
following residue formula (Theorem 4.4 of \cite{S}):
$$
\sum_{n\in N_{reg}} e^{\langle t,2i\pi n\rangle} f(2i\pi n)
=\sum_{o\sigma\in OB} Res^{o\sigma}
\left( e^{\langle t,z\rangle}f(z)F^{\mathfrak a}_{\sigma}(-z)\right).
$$
Thus, when 
$$
f=\frac{1}{\prod_{j=1}^k \alpha_j}
$$ 
is sufficiently decreasing, this formula expresses the series 
$$
\sum_{n\in \Z^r,\langle \alpha_j,n\rangle\neq 0} 
\frac{1}{\prod_{j=1}^k \langle \alpha_j,2i\pi n\rangle}
$$ 
as an explicit rational number.  
\bigskip

\begin{proof}
From the definitions of $Z(f)(t,z)$ and $CT(f)(t)$,
we obtain for any $P\in S(V^*)$:
$$
P(\partial_t)Z(f)(t,0)=Z(Pf)(t,0),~ 
P(\partial_t)CT(f)(t)=CT(Pf)(t).
$$
Thus, it is enough to prove that 
$Z(f)(t,0)=CT(f)(t)$ for $f\in G_{\Delta}$, because $G_{\Delta}$
generates $R_{\Delta}$ as a $S(V^*)$-module by Lemma \ref{mod}.

For $t$ in an alcove ${\mathfrak a}$, we can define the operator
$Z^t:R_{\Delta}\to \CO$ by
$$
Z^t(f)(z)=\sum_{n\in N_{reg}}e^{\langle t,z+2i\pi n\rangle}f(z+2i\pi n).
$$
The kernel formula holds for the operator $Z^t$. In particular, we obtain 
for $f\in G_{\Delta}$:
$$
Z^t(f)(0)=Tr_{S_\Delta}\left(Res_{\Delta}m(e^{\langle t,\cdot\rangle}f)
C Z^t Res_{\Delta}\right).
$$
We thus need to prove that, for $f\in G_{\Delta}$,
$$
Tr_{S_\Delta}\left( Res_{\Delta}m(e^{\langle t,\cdot\rangle }f)
C(E^t- Z^t) Res_{\Delta}\right)=0.
$$
But $E^t$ is given by a sum over the full lattice $N$, while $Z^t$ is
only over the regular elements of $N$. Thus, we can write (in many
ways) $E^t-Z^t$ as a linear combination of  
operators $E^t_{\Gamma_\alpha}$ where each $\Gamma_\alpha$ is an
admissible subset of $N$ contained in the real hyperplane $H_{\alpha}$.
Now Szenes formula will follow from
\begin{proposition}
Let $\Gamma$ be an admissible subset of $N$ contained in the real
hyperplane $H_{\alpha}$. Then, for $f\in G_{\Delta}$ 
$$
Tr_{S_\Delta}\left( Res_{\Delta}m( e^{\langle t,\cdot\rangle }f)
C E_\Gamma^t Res_{\Delta}\right)=0.
$$ 
\end{proposition}
\begin{proof}
It suffices to prove that 
$$
\sum_{o\sigma\in OB} Res^{o\sigma}\left(e^{\langle t,z\rangle}
f(z)E_{\Gamma}^t(\phi_{\sigma})(-z)\right) = 0
$$
for some diagonal basis $OB$.

A total order on $\Delta$ provides us with a special diagonal basis
$OB$ of $O\CB(\Delta)$ (see for example \cite{BV}, Proposition 14.) We
choose this order 
such that $\alpha$ is minimal. In this case every element of $OB$ is
of the form $o\sigma=(\alpha_1,\alpha_2,\ldots,\alpha_r)$ with
$\alpha_1=\alpha$. We claim that for each $o\sigma\in OB$,
$$
Res^{o\sigma}\left(e^{\langle t,z\rangle}f(z)
E_{\Gamma}^t(\phi_{\sigma})(-z)\right)=0.
$$
Indeed, we use the notation of Lemma \ref{b} and write
$V'=H_{\alpha}$. Then our set $\Gamma$ is contained in $V'$. Thus
$$
E_{\Gamma}^t(\phi_{\sigma})(z_1,z')=\frac{e^{t_1z_1}}{z_1}
\sum_{\gamma\in \Gamma}\frac{e^{\langle t',z'+2i\pi\gamma\rangle}}
{\prod_{j=2}^r\langle \alpha_j,z'+2i\pi \gamma\rangle}.
$$
We see that for $t$ fixed and regular,
$$
e^{\langle t,z\rangle }f(z)E_{\Gamma}^t(\phi_{\sigma})(-z)=
\frac{1}{z_1}f(z_1,z')\psi(z')
$$
where $f\in G_{\Delta}$ and $\psi(z')$ has poles at most on the
complex hyperplanes $\alpha_j=0$ for $j=2,\cdots,r$. Thus the claim
follows from Lemma \ref{b}.
\end{proof}
\end{proof}

\end{document}